\documentclass{amsart}
\usepackage{amsmath,amsfonts,amssymb,amsthm}
\input amssym.def
\def \Z{\Bbb Z}
\def \C{\Bbb C}

\def \N{\Bbb N}

\def \D{{\mathcal{D}}}
\def \E{{\mathcal{E}}}
\def \E{{\mathcal{E}}}

\def \Res{{\rm Res}}

\def \End{{\rm End}}

\def \Aut{{\rm Aut}}

\def \Hom{{\rm Hom}}

\def \<{\langle}
\def \>{\rangle}

\def \g{{\frak{g}}}

\def \gl{{\frak{gl}}}

\def \be{\begin{equation}\label}
\def \ee{\end{equation}}
\def \bex{\begin{example}\label}
\def \eex{\end{example}}
\def \bl{\begin{lem}\label}
\def \el{\end{lem}}
\def \bt{\begin{thm}\label}
\def \et{\end{thm}}
\def \bp{\begin{prop}\label}
\def \ep{\end{prop}}
\def \br{\begin{rem}\label}
\def \er{\end{rem}}
\def \bc{\begin{coro}\label}
\def \ec{\end{coro}}
\def \bd{\begin{de}\label}
\def \ed{\end{de}}

\newtheorem{thm}{Theorem}[section]
\newtheorem{prop}[thm]{Proposition}
\newtheorem{coro}[thm]{Corollary}

\newtheorem{example}[thm]{Example}
\newtheorem{lem}[thm]{Lemma}
\newtheorem{rem}[thm]{Remark}
\newtheorem{de}[thm]{Definition}

\makeatletter \@addtoreset{equation}{section}

\makeatother \makeatletter

\begin{document}
\title{On quasi modules at infinity for vertex algebras}

\author{Haisheng Li}
\address{HL: Department of Mathematical Sciences, Rutgers University, Camden,
NJ 08102, USA}
\email{hli@camden.rutgers.edu}
\thanks{*Corresponding author}

\author{Qiang Mu*}
\address{QM: School of Mathematical Sciences, Harbin Normal University,
Harbin, Heilongjiang 150080, China}
\email{qmu520@gmail.com}





\maketitle

\begin{abstract}
A theory of quasi modules at infinity for (weak) quantum vertex algebras including vertex algebras
was previously developed in \cite{li-infinity}.
In this current paper, quasi modules at infinity for vertex algebras are revisited.
Among the main results, we extend some technical results, to fill in a gap in the proof of a theorem therein,
and we obtain a commutator formula for general quasi modules at infinity
and establish a version of the converse of the aforementioned theorem.
\end{abstract}

\section{Introduction}
In order to associate quantum vertex algebras to certain algebras such as centerless double Yangians
and the Lie algebra of pseudo-differential operators on the circle,
a theory of what were called quasi modules at infinity for quantum vertex algebras was developed
in \cite{li-infinity}. The notion of (weak) quantum vertex algebra, which was formulated in \cite{li-qva1} (cf. \cite{ek}),
is a natural generalization of the notions of vertex algebra
and vertex super-algebra.
Let $V$ be a weak quantum vertex algebra. Note that for a (quasi) $V$-module $W$,
each vector $v\in V$ is represented by a formal series $Y_{W}(v,x)\in \Hom (W,W((x)))$.
In contrast, for a (quasi) $V$-module at infinity $M$ each vector $v\in V$
is represented by a formal series $Y_{M}(v,x)\in \Hom (M,M((x^{-1})))$.
Normally, associative or Lie algebras are associated with vertex algebras or more general quantum vertex algebras through
their ``highest weight'' modules. A matter of fact is that certain algebras
such as the Lie algebra of pseudo-differential operators on the circle only admit  ``lowest weight'' modules;
they do {\em not} admit nontrivial highest weight modules. This was the main motivation
for developing the theory of quasi modules at infinity. Indeed, this theory enabled us
to associate quantum vertex algebras to the aforementioned algebras
through their lowest weight modules.

Let $V$ be a weak quantum vertex algebra.
A {\em $V$-module at infinity} is a vector space $W$ equipped with a linear map
\begin{eqnarray*}
Y_{W}(\cdot,x):&& V\rightarrow \Hom (W,W((x^{-1})))\\
&&v\mapsto Y_{W}(v,x),
\end{eqnarray*}
satisfying the conditions that $Y_{W}({\bf 1},x)=1_{W}$ and that for any $u,v\in V$,
there exists a nonnegative integer $k$ such that
$$(x_{1}-x_{2})^{k}Y_{W}(u,x_{1})Y_{W}(v,x_{2})\in \Hom \left(W,W((x_{1}^{-1},x_{2}^{-1}))\right),$$
\begin{eqnarray}
\mbox{}\ \ \ \ \ \ \ x_{0}^{k}Y_{W}(Y(u,x_{0})v,x_{2})=\left((x_{1}-x_{2})^{k}Y_{W}(u,x_{1})Y_{W}(v,x_{2})\right)|_{x_{1}=x_{2}+x_{0}}.
\end{eqnarray}
For the definition of a quasi $V$-module at infinity, the second condition above is replaced with
the condition that for any $u,v\in V$,
there is a nonzero polynomial $p(x_{1},x_{2})$ such that
$$p(x_{1},x_{2})Y_{W}(u,x_{1})Y_{W}(v,x_{2})\in \Hom \left(W,W((x_{1}^{-1},x_{2}^{-1}))\right),$$
\begin{eqnarray*}
p(x_{0}+x_{2},x_{2})Y_{W}(Y(u,x_{0})v,x_{2})=\left(p(x_{1},x_{2})Y_{W}(u,x_{1})Y_{W}(v,x_{2})\right)|_{x_{1}=x_{2}+x_{0}}.
\end{eqnarray*}
It was proved therein that if $V$ is a vertex algebra, the second condition in the definition of a module at infinity
amounts to the following opposite Jacobi identity
\begin{eqnarray*}
&&x_{0}^{-1}\delta\left(\frac{x_{1}-x_{2}}{x_{0}}\right)Y_{W}(v,x_{2})Y_{W}(u,x_{1})
-x_{0}^{-1}\delta\left(\frac{x_{2}-x_{1}}{-x_{0}}\right)Y_{W}(u,x_{1})Y_{W}(v,x_{2})\nonumber\\
&&\hspace{1cm}=x_{2}^{-1}\delta\left(\frac{x_{1}-x_{0}}{x_{2}}\right)Y_{W}(Y(u,x_{0})v,x_{2}).
\end{eqnarray*}
That is, for a vertex algebra $V$, the notion of
$V$-module at infinity is the same as that of right $V$-module
(see \cite{hl}, \cite{li-reg}, \cite{xue}).

A notion of quasi module for vertex algebras was introduced in \cite{li-new}, in order to associate
vertex algebras to certain Lie algebras. Meanwhile,
a notion of vertex $\Gamma$-algebra was also introduced with $\Gamma$ a group, where
a {\em vertex $\Gamma$-algebra} is a vertex algebra $V$ equipped with two group homomorphisms
$$L:\ \Gamma \rightarrow GL(V)\ \mbox{ and }\ \phi:\ \Gamma \rightarrow \C^{\times},$$
satisfying the condition that $L(g){\bf 1}={\bf 1}$,
$$L(g)Y(v,x)=Y(L(g)v,\phi(g)x)L(g)\ \ \ \mbox{ for }g\in \Gamma, \ v\in V.$$
 In \cite{li-infinity}, a notion of quasi $V$-module at infinity for a vertex $\Gamma$-algebra $V$ was studied.
For simplicity, consider the special case with $\Gamma$ a group of linear functions
$g(x)=\alpha x+\beta$ with $\alpha\in \C^{\times}, \beta\in \C$ (with respect to function composition).
In this case, a quasi $V$-module at infinity is a quasi module at infinity $(W,Y_{W})$ for $V$ viewed as a vertex algebra,
which also satisfies
$$Y_{W}(L(g)v,x)=Y_{W}(v,g^{-1}(x))\ \ \ \mbox{ for }g\in \Gamma,\ v\in V.$$

In this current paper, we revisit quasi modules at infinity for vertex algebras.
As one of the main purposes of this paper, we generalize a technical result of \cite{li-infinity} (Lemma 5.11),
to fill in a gap in the proof of Theorem 5.14 therein.
(Roughly speaking, this theorem asserts that each restricted module for a certain Lie algebra
with some parameters is naturally a quasi module at infinity for a certain vertex $\Gamma$-algebra.)
As the main results of this paper, for a general vertex $\Gamma$-algebra $V$ we obtain a commutator formula
for quasi $V$-modules at infinity, which to a certain extent is analogous to the twisted vertex operator commutator formula
(see \cite{flm}). (Note that for {\em quasi} modules at infinity we {\em no longer} have the opposite Jacobi identity.)
Just as the twisted commutator formula is important in the study of twisted modules,
this commutator formula is important in the study of quasi modules at infinity for vertex $\Gamma$-algebras.
As an application, we establish a version of the converse of the aforementioned theorem.

This paper is organized as follows: In Section 2, we refine or extend several results of \cite{li-infinity}
and we present a commutator formula.
In Section 3, we present a complete proof of Theorem 5.14 of \cite{li-infinity} and
establish a converse.

\section{Quasi modules at infinity for vertex algebras}
In this section, we first recall some basic results from \cite{li-infinity} and
we then refine several results and establish a commutator formula for quasi modules at infinity.

We begin with some basic formal variable notations. First of all, throughout this paper,
vector spaces will be over $\C$ (the field of complex numbers) unless it is stated otherwise.
We shall use the formal variable notations and conventions as established in \cite{flm} and \cite{fhl}.
In particular, for a vector space $W$,
$W[[x,x^{-1}]]$ and $W[[x_{1}^{\pm 1},x_{2}^{\pm 1}]]$ denote the spaces of doubly infinite formal series with coefficients in $W$,
while
$W((x))$ and $W((x_{1},x_{2}))$ denote the spaces of lower truncated infinite formal series.
Recall also the formal delta functions:
$$\delta(x)=\sum_{n\in \Z}x^{n},$$
$$x_{1}^{-1}\delta\left(\frac{x_{2}}{x_{1}}\right)=\sum_{n\in \Z}x_{1}^{-n-1}x_{2}^{n}.$$

Denote by $\C(x)$ the field of rational functions in variable $x$
and by $\C(x_{1},x_{2})$ the field of rational functions in $x_{1}$ and $x_{2}$.
We have canonical field embeddings:
\begin{eqnarray*}
\iota_{x_{1},x_{2}}:&&  \C(x_{1},x_{2})\rightarrow \C((x_{1}))((x_{2})),\\
\iota_{x_{1}^{-1},x_{2}}:&& \C(x_{1},x_{2})\rightarrow \C((x_{1}^{-1}))((x_{2})),
\end{eqnarray*}
which extend the ring embeddings of $\C[x_{1},x_{2}]$ into $\C((x_{1}))((x_{2}))$
and $\C((x_{1}^{-1}))((x_{2}))$, respectively.

Throughout this paper, we set
\begin{eqnarray}
G=\{ \alpha x+\beta\ |\ \alpha,\beta\in \C,\ \alpha\ne 0\}\subset \C[x],
\end{eqnarray}
which is considered as a group with respect to function composition.
(Note that $G$ consists of all the linear transformations that preserve $\infty$.)
For $g(x)=\alpha x+\beta\in G$, we have $g^{-1}(x)=\alpha^{-1}x-\alpha^{-1}\beta$.
We have a canonical group homomorphism
\begin{eqnarray}
\Phi:\  G\rightarrow \C^{\times}, \ \ g(x)=\alpha x+\beta\mapsto \alpha.
\end{eqnarray}

Let $g(x)=\alpha x+\beta\in G$. For $n\in \Z$, we view $g(x)^{n}$ as an element of $\C((x^{-1}))$:
$$g(x)^{n}=(\alpha x+\beta)^{n}=\sum_{i\ge 0}\binom{n}{i}\alpha^{n-i}\beta^{i}x^{n-i}.$$
It is understood that
\begin{eqnarray}
\delta\left(\frac{g(x_{2})}{x_{1}}\right)=\sum_{n\in \Z}x_{1}^{-n}g(x_{2})^{n}=\sum_{n\in \Z}\sum_{i\ge 0}\binom{n}{i}\alpha^{n-i}\beta^{i}x_{1}^{-n}x_{2}^{n-i}.
\end{eqnarray}

For a vector space $W$ over $\C$, set
\begin{eqnarray}
\E^{o}(W)=\Hom (W,W((x^{-1})))\subset (\End W)[[x,x^{-1}]].
\end{eqnarray}
Vector space $\E^{o}(W)$ (over $\C$) is also naturally a module over $\C((x^{-1}))$.
The identity operator, denoted by $1_{W}$, is a special element of $\E^{o}(W)$.

Let $a(x)=\sum_{n\in \Z}a_{n}x^{-n-1}\in \E^{o}(W)$ and let $g(x)=\alpha x+\beta\in G$.
As a convention we define
\begin{eqnarray*}
\mbox{}\ \ \ \ \  a(g(x))=\sum_{n\in \Z} a_{n}(\alpha x+\beta)^{-n-1}
=\sum_{n\in \Z}\sum_{i\ge 0}\binom{-n-1}{i}\alpha^{-n-i-1}\beta^{i}a_{n}x^{-n-i-1},
\end{eqnarray*}
which still lies in $\E^{o}(W)$.

For $g\in G$, define $L_{g}\in \End (\E^{o}(W))$ by
\begin{eqnarray}
L_{g}(a(x))=a(g^{-1}(x))\ \ \ \mbox{ for }g\in G,\ a(x)\in \E^{o}(W).
\end{eqnarray}
This makes $\E^{o}(W)$ a (left) $G$-module\footnote{An action of $G$ on $\E^{o}(W)$
was defined in \cite{li-infinity} by $R_{g}(a(x))=a(g(x))$ for $g\in G,\ a(x)\in \E^{o}(W)$.
In fact, this gives a {\em right} action instead of a {\em left} action because $G$ is nonabelian.}.

\bd{dquasi-module-infinity}
{\em Let $V$ be a vertex algebra. A {\em quasi $V$-module at infinity} is a vector space $W$
equipped with a linear map
\begin{eqnarray*}
Y_{W}(\cdot,x):&& V\rightarrow \Hom (W,W((x^{-1})))\subset (\End W)[[x,x^{-1}]]\\
&&v\mapsto Y_{W}(v,x)=\sum_{n\in \Z}v_{n}x^{-n-1}\ \ (\mbox{with }v_{n}\in \End W),
\end{eqnarray*}
satisfying the conditions that $Y_{W}({\bf 1},x)=1_{W}$ (the identity operator on $W$) and
that for any $u,v\in V$, there exists a nonzero polynomial $p(x_{1},x_{2})$
such that
\begin{eqnarray}\label{econdition-puv}
p(x_{1},x_{2})Y_{W}(u,x_{1})Y_{W}(v,x_{2})\in \Hom \left(W,W((x_{1}^{-1},x_{2}^{-1}))\right)
\end{eqnarray}
and
\begin{eqnarray}\label{eweak-assoc-def}
&&\left(p(x_{1},x_{2})Y_{W}(u,x_{1})Y_{W}(v,x_{2})\right)|_{x_{1}=x_{2}+x_{0}}=
p(x_{0}+x_{2},x_{2})Y_{W}(Y(u,x_{0})v,x_{2}).\nonumber\\
&&
\end{eqnarray}}
\ed

\br{rafter-definition} {\em Note that for $A(x_{1},x_{2})\in \Hom (W,W((x_{1}^{-1},x_{2}^{-1})))$ with $W$ a vector space,
$$A(x_{2}+x_{0},x_{2})\ \ \mbox{ exists in }\left(\Hom (W,W((x_{2}^{-1})))\right)[[x_{0}]].$$
The condition (\ref{econdition-puv}) in Definition \ref{dquasi-module-infinity}
guarantees that the expression on the left-hand side of (\ref{eweak-assoc-def}) exists.
On the other hand, assume that $(W,Y_{W})$ is a quasi module at infinity for a vertex algebra $V$ and
let $u,v\in V$. Then (\ref{eweak-assoc-def}) holds for any  nonzero polynomial $p(x_{1},x_{2})$
such that (\ref{econdition-puv}) holds.}
\er

We have the following result which strengthens Lemma 5.1 of \cite{li-infinity}:

\bp{pjacobi} Let $(W,Y_{W})$ be a quasi module at infinity for a vertex algebra $V$.
Then for any $u,v\in V$, there exists a nonzero polynomial $p(x_{1},x_{2})$
such that
\begin{eqnarray}\label{epquasi-local}
p(x_{1},x_{2})Y_{W}(v,x_{2})Y_{W}(u,x_{1})=p(x_{1},x_{2})Y_{W}(u,x_{1})Y_{W}(v,x_{2})
\end{eqnarray}
and
\begin{eqnarray}\label{ep-jacobi}
&&x_{0}^{-1}\delta\left(\frac{x_{1}-x_{2}}{x_{0}}\right)p(x_{1},x_{2})Y_{W}(v,x_{2})Y_{W}(u,x_{1})\\
&&\hspace{1cm}-x_{0}^{-1}\delta\left(\frac{x_{2}-x_{1}}{-x_{0}}\right)p(x_{1},x_{2})Y_{W}(u,x_{1})Y_{W}(v,x_{2})
\nonumber\\
&=&x_{2}^{-1}\delta\left(\frac{x_{1}-x_{0}}{x_{2}}\right)p(x_{1},x_{2})Y_{W}(Y(u,x_{0})v,x_{2}).\nonumber
\end{eqnarray}
\ep

\begin{proof} Lemma 5.1 in \cite{li-infinity} asserts that
there exists a nonzero polynomial $p(x_{1},x_{2})$ such that (\ref{epquasi-local}) holds, from which
we see that (\ref{econdition-puv}) holds.
{}From Remark \ref{rafter-definition}, (\ref{eweak-assoc-def}) holds.
Then (\ref{ep-jacobi}) follows immediately from Lemma 4.2 therein.
\end{proof}

\bd{dquasi-local} {\em Let $W$ be a vector space.
A subset $U$ of $\E^{o}(W)$ is said to be {\em quasi local} if for any $a(x),b(x)\in U$,
there exists a nonzero polynomial $p(x_{1},x_{2})$ such that
\begin{eqnarray}\label{edef-pab=0}
p(x_{1},x_{2})[a(x_{1}),b(x_{2})]=0.
\end{eqnarray}}
\ed

The following notion was due to \cite{gkk}:

\bd{dgamma-local} {\em Let $W$ be a vector space and let $\Gamma$ be a subgroup of $G$.
A subset $U$ of $\E^{o}(W)$ is said to be {\em $\Gamma$-local} if for any $a(x),b(x)\in U$,
(\ref{edef-pab=0}) holds with
$p(x_{1},x_{2})$ a product of $x_{1}-g(x_{2})$ with $g\in \Gamma$ (not necessarily multiplicity-free).}
\ed

Note that the relation (\ref{edef-pab=0}), which can be written as
$$p(x_{1},x_{2})a(x_{1})b(x_{2})=p(x_{1},x_{2})b(x_{2})a(x_{1}),$$
 implies
\begin{eqnarray}\label{epab-condition}
p(x_{1},x_{2})a(x_{1})b(x_{2})\in \Hom \left(W,W((x_{1}^{-1},x_{2}^{-1}))\right).
\end{eqnarray}

A subset $U$ of $\E^{o}(W)$ is said to be {\em quasi compatible} if
 for any $a(x),b(x)\in U$, there exists a nonzero polynomial
$p(x_{1},x_{2})$ such that (\ref{epab-condition}) holds. Furthermore,
for a subgroup $\Gamma$ of $G$ we define {\em $\Gamma$-compatibility} correspondingly.

Now, let $a(x),b(x)\in \E^{o}(W)$. Assume that there exists a nonzero polynomial
$p(x_{1},x_{2})$ such that (\ref{epab-condition}) holds.
We define
$$a(x)_{n}b(x)\in \E^{o}(W)\ \ \ \mbox{ for }n\in \Z$$
in terms of generating function
$$Y_{\E^{o}}(a(x),x_{0})b(x)=\sum_{n\in \Z}a(x)_{n}b(x)x_{0}^{-n-1}$$
by
\begin{eqnarray*}
Y_{\E^{o}}(a(x),x_{0})b(x)=\iota_{x^{-1}, x_{0}}\left(\frac{1}{p(x+x_{0},x)}\right)
\left(p(x_{1},x)a(x_{1})b(x)\right)|_{x_{1}=x+x_{0}},
\end{eqnarray*}
where $\iota_{x^{-1}, x_{0}}:  \C(x,x_{0})\rightarrow \C((x^{-1}))((x_{0}))$
is the unique extension of the embedding of $\C[x,x_{0}]$ into $\C((x^{-1}))((x_{0}))$.

The following result was obtained in \cite{li-infinity} (Theorem 5.4):

\bt{trecall-main}
Let $W$ be a vector space and let $U$ be any quasi local subset of $\E^{o}(W)$.
Then $U$ generates a vertex algebra $\<U\>$ with $Y_{\E^{o}}$ as the vertex operator map and
with $1_{W}$ as the vacuum vector, and $W$ is a quasi module at infinity with
$$Y_{W}(\alpha(x),z)=\alpha(z)\ \ \ \mbox{ for }\alpha(x)\in \<U\>.$$
\et

The following was also essentially obtained in \cite{li-infinity}:

\bp{prepresentation}
Let $V$ be a vertex algebra and let $W$ be a vector space equipped with a linear map
$Y_{W}(\cdot,x): V\rightarrow \Hom (W,W((x^{-1})))$ with $Y_{W}({\bf 1},x)=1_{W}$. Set
$$V_{W}=\{ Y_{W}(v,x)\ |\ v\in V\}\subset \E^{o}(W).$$
Then
$(W,Y_{W})$ is a quasi $V$-module at infinity if and only if $V_{W}$ is quasi local,
$(V_{W},Y_{\E^{o}},1_{W})$ carries the structure of a vertex algebra,
 and $Y_{W}$ is a vertex algebra homomorphism from $V$ to $V_{W}$.
\ep

\begin{proof} Assume that $(W,Y_{W})$ is a quasi $V$-module at infinity.
By Proposition \ref{pjacobi},
$V_{W}$ is a quasi local subspace of $\E^{o}(W)$.
{}From \cite{li-infinity} (the first half of the proof of Lemma 5.9) we have
\begin{eqnarray}\label{ehom-property-1}
Y_{\E^{o}}(Y_{W}(u,x),x_{0})Y_{W}(v,x)=Y_{W}(Y(u,x_{0})v,x)\ \ \ \mbox{ for }u,v\in V.
\end{eqnarray}
We also have $Y_{W}({\bf 1},x)=1_{W}$.
Then it follows that $(V_{W},Y_{\E^{o}},1_{W})$ carries the structure of a vertex algebra and $Y_{W}$ is a homomorphism of vertex algebras.
The other direction is clear from Theorem \ref{trecall-main}.
\end{proof}

The following result generalizes Lemma 5.11 of \cite{li-infinity}:

\bl{lextension}
Let $W$ be a vector space and let $a(x),b(x)\in \E^{o}(W)$. Assume
\begin{eqnarray}
&&[a(x_{1}),b(x_{2})]=\sum_{j=0}^{r}A_{j}(x_{2})\frac{1}{j!}\left(\frac{\partial}{\partial x_{2}}\right)^{j}x_{1}^{-1}\delta\left(\frac{x_{2}}{x_{1}}\right)\\
&&\hspace{2cm}+\sum_{i=1}^{k}\sum_{j=0}^{s}B_{ij}(x_{2})\frac{1}{j!}\left(\frac{\partial}{\partial x_{2}}\right)^{j}x_{1}^{-1}\delta\left(\frac{g_{i}(x_{2})}{x_{1}}\right),\nonumber
\end{eqnarray}
where $A_{j}(x), B_{ij}(x)\in \E^{o}(W)$, $g_{i}(x)\in G$ with $g_{i}(x)\ne x$ for $1\le i\le k$. Then
\begin{eqnarray}
&&a(x)_{j}b(x)=-A_{j}(x)\ \ \ \mbox{ for }0\le j\le r,\nonumber\\
&&a(x)_{j}b(x)=0\ \ \ \mbox{ for }j>r.
\end{eqnarray}
\el

\begin{proof}  Set
$$q(x_{1},x)=\prod_{i=1}^{k}\left((x_{1}-x)^{r+1}-(g_{i}(x)-x)^{r+1}\right)^{s+1}\in \C[x_{1},x].$$
Noticing that
$$q(x_{1},x)=\bar{q}(x_{1},x)\prod_{i=1}^{k}(x_{1}-g_{i}(x))^{s+1}$$
for some $\bar{q}(x_{1},x)\in \C[x_{1},x]$ and that
$$(x_{1}-g_{i}(x))^{s+1}\left(\frac{\partial}{\partial x_{2}}\right)^{j}x_{1}^{-1}\delta\left(\frac{g_{i}(x_{2})}{x_{1}}\right)=0$$
for $1\le i\le k,\ 0\le j\le s$, we have
$$q(x_{1},x_{2})\left(\frac{\partial}{\partial x_{2}}\right)^{j}x_{1}^{-1}\delta\left(\frac{g_{i}(x_{2})}{x_{1}}\right)=0.$$
Then we get
\begin{eqnarray}\label{eqab=}
\mbox{}\ \ \ \ \ \ \ q(x_{1},x_{2})[a(x_{1}),b(x_{2})]=\sum_{j=0}^{r}q(x_{1},x_{2})A_{j}(x_{2})\frac{1}{j!}\left(\frac{\partial}{\partial x_{2}}\right)^{j}x_{1}^{-1}\delta\left(\frac{x_{2}}{x_{1}}\right).
\end{eqnarray}

Note that (\ref{eqab=}) implies
$$(x_{1}-x_{2})^{r+1}q(x_{1},x_{2})[a(x_{1}),b(x_{2})]=0.$$
{}From the definition of $Y_{\E^{o}}$ we have
\begin{eqnarray*}
x_{0}^{r+1}q(x+x_{0},x)Y_{\E^{o}}(a(x),x_{0})b(x)=\left((x_{1}-x)^{r+1}q(x_{1},x)a(x_{1})b(x)\right)|_{x_{1}=x+x_{0}}.
\end{eqnarray*}
Combining the two identities above with Lemma 4.2 of \cite{li-infinity},  we obtain
\begin{eqnarray*}
&&x_{0}^{-1}\delta\left(\frac{x_{1}-x}{x_{0}}\right)(x_{1}-x)^{r+1}q(x_{1},x)b(x)a(x_{1})\hspace{2cm}\nonumber\\
&&\ \ \ \ -x_{0}^{-1}\delta\left(\frac{x-x_{1}}{-x_{0}}\right)(x_{1}-x)^{r+1}q(x_{1},x)a(x_{1})b(x)\nonumber\\
&=&x^{-1}\delta\left(\frac{x_{1}-x_{0}}{x}\right)(x_{1}-x)^{r+1}q(x_{1},x)Y_{\E^{o}}(a(x),x_{0})b(x).
\end{eqnarray*}
Applying $\Res_{x_{0}}x_{0}^{-r-1}$ to both sides, we get
\begin{eqnarray*}
-q(x_{1},x)[a(x_{1}),b(x)]=\sum_{i\ge 0}a(x)_{i}b(x) \frac{1}{i!}q(x_{1},x)\left(\frac{\partial}{\partial x}\right)^{i}x_{1}^{-1}\delta\left(\frac{x}{x_{1}}\right),
\end{eqnarray*}
which is a finite sum.
Combining this with (\ref{eqab=}) we obtain
\begin{eqnarray}\label{efinalcomp}
\sum_{i\ge 0}\left(A_{i}(x)+a(x)_{i}b(x)\right)\frac{1}{i!}q(x_{1},x)\left(\frac{\partial}{\partial x}\right)^{i}x_{1}^{-1}\delta\left(\frac{x}{x_{1}}\right)=0,
\end{eqnarray}
where we set $A_{i}(x)=0$ for $i>r$.
Notice that we have
\begin{eqnarray*}
q(x_{1},x)=(x_{1}-x)^{r+1}P(x_{1},x)+Q(x),
\end{eqnarray*}
where $P(x_{1},x)$ is a polynomial and
$$Q(x)=(-1)^{k(s+1)}\prod_{i=1}^{k}(g_{i}(x)-x)^{(r+1)(s+1)}.$$
Since
$$(x_{1}-x)^{r+1}\left(\frac{\partial}{\partial x}\right)^{j}x_{1}^{-1}\delta\left(\frac{x}{x_{1}}\right)=0
\ \ \ \mbox{ for }0\le j\le r,$$
 (\ref{efinalcomp}) reduces to
\begin{eqnarray*}\label{elastcomp}
\sum_{i\ge 0}\left(A_{i}(x)+a(x)_{i}b(x)\right)\frac{1}{i!}Q(x)\left(\frac{\partial}{\partial x}\right)^{i}x_{1}^{-1}\delta\left(\frac{x}{x_{1}}\right)=0.
\end{eqnarray*}
{}From \cite{li-local} (Lemma 2.1.4), we get
$$\left(A_{i}(x)+a(x)_{i}b(x)\right)Q(x)=0\ \ \ \mbox{ for }i\ge 0.$$
Since $g_{i}(x)\ne x$ for $1\le i\le k$, we have $Q(x)\ne 0$. Then our assertions follow immediately.
\end{proof}

As a generalization and a corollary of Lemma \ref{lextension} we have:

\bp{pextension}
Let $W$ be a vector space and let $a(x),b(x)\in \E^{o}(W)$. Assume
\begin{eqnarray}\label{eab=given}
[a(x_{1}),b(x_{2})]=
\sum_{i=1}^{k}\sum_{j=0}^{r}A_{ij}(x_{2})\frac{1}{j!}\left(\frac{\partial}{\partial x_{2}}\right)^{j}x_{1}^{-1}\delta\left(\frac{g_{i}(x_{2})}{x_{1}}\right),
\end{eqnarray}
where $A_{ij}(x)\in \E^{o}(W)$ and $g_{i}(x)\in G$ distinct for $1\le i\le k$. Then
\begin{eqnarray}
&&\Phi(g_{i})a(g_{i}(x))_{j}b(x)=-A_{ij}(x)\ \ \ \mbox{ for }0\le j\le r,\nonumber\\
&&\Phi(g_{i})a(g_{i}(x))_{j}b(x)=0\ \ \ \mbox{ for }j>r,
\end{eqnarray}
where
$\Phi: G\rightarrow \C^{\times}$ was defined before by $\Phi(g)=\alpha$ for $g(x)=\alpha x+\beta \in G$.
\ep

\begin{proof} Let $g(x)=\alpha x+\beta\in G$ (with $\alpha\in \C^{\times},\ \beta\in \C$). We have
\begin{eqnarray*}
&&x_{1}^{-1}\delta\left(\frac{g(x_{2})}{x_{1}}\right)=x_{1}^{-1}\delta\left(\frac{\alpha x_{2}+\beta}{x_{1}}\right)
=\alpha^{-1}x_{2}^{-1}\delta\left(\frac{x_{1}-\beta}{\alpha x_{2}}\right)\\
&=&\alpha^{-1}x_{2}^{-1}\delta\left(\frac{\alpha^{-1}x_{1}-\alpha^{-1}\beta}{x_{2}}\right)
=\Phi(g)^{-1}x_{2}^{-1}\delta\left(\frac{g^{-1}(x_{1})}{x_{2}}\right).
\end{eqnarray*}
Using this, {}from (\ref{eab=given}) we have
\begin{eqnarray*}
[a(g(x_{1})),b(x_{2})]=
\sum_{i=1}^{k}\sum_{j=0}^{r}\Phi(g_{i})^{-1}A_{ij}(x_{2})\frac{1}{j!}\left(\frac{\partial}{\partial x_{2}}\right)^{j}x_{2}^{-1}\delta\left(\frac{g_{i}^{-1}g(x_{1})}{x_{2}}\right).
\end{eqnarray*}
Taking $g(x)=g_{i}(x)$ with $1\le i\le k$, by Lemma \ref{lextension} we obtain
$$a(g_{i}(x))_{j}b(x)=-\Phi(g_{i})^{-1}A_{ij}(x)\ \ \ \mbox{ for } 0\le j\le r$$
and $a(g_{i}(x))_{j}b(x)=0$ for $j>r$,
as desired.
\end{proof}

The following technical result generalizes a result of \cite{gkk}:

\bl{lformal-variable}
Let $W$ be a vector space, let
$$A(x_{1},x_{2})\in \Hom \left(W,W((x_{2}^{-1}))[[x_{1},x_{1}^{-1}]]\right),$$
and let
$$p(x_{1},x_{2})=\prod_{i=1}^{r}(x_{1}-g_{i}(x_{2}))^{k_{i}},$$
where $g_{1}(x),\dots,g_{r}(x)$ are distinct elements of $G$ and
$k_{1},\dots,k_{r}$ are positive integers. Then
$p(x_{1},x_{2})A(x_{1},x_{2})=0$ if and only if
\begin{eqnarray}
A(x_{1},x_{2})=\sum_{i=1}^{r}\sum_{j=0}^{k_{i}-1}A_{ij}(x_{2})\left(\frac{\partial}{\partial x_{2}}\right)^{j}x_{1}^{-1}\delta\left(\frac{g_{i}(x_{2})}{x_{1}}\right)
\end{eqnarray}
for some $A_{ij}(x)\in \E^{o}(W)$.
\el

\begin{proof}  The ``if'' part is clear. For the ``only if'' part, we first consider the case $r=1$.
Note that for any $h(x)\in G$, $A(x_{1},h(x_{2}))$ exists in $\Hom \left(W,W((x_{2}^{-1}))[[x_{1},x_{1}^{-1}]]\right)$.
Let us simply use $g(x)$ for $g_{1}(x)$.  We see that
$$(x_{1}-g(x_{2}))^{k}A(x_{1},x_{2})=0$$
 if and only if
$$(x_{1}-x_{2})^{k}A(x_{1},g^{-1}(x_{2}))=0.$$
It was known (see \cite{kac2}, \cite{dlm}) that the latter is equivalent to
\begin{eqnarray}\label{eagb-exp}
A(x_{1},g^{-1}(x_{2}))=\sum_{j=0}^{k-1} B_{j}(x_{2})
\left(\frac{\partial}{\partial x_{2}}\right)^{j}x_{1}^{-1}\delta\left(\frac{x_{2}}{x_{1}}\right)
\end{eqnarray}
for some $B_{j}(x)\in (\End W)[[x,x^{-1}]]$ . For $0\le j\le k-1$, since
$$\Res_{x_{1}}(x_{1}-x_{2})^{j}\left(\frac{\partial}{\partial x_{2}}\right)^{j}x_{1}^{-1}\delta\left(\frac{x_{2}}{x_{1}}\right)=j!,$$
we have
$$B_{j}(x_{2})=\frac{1}{j!}\Res_{x_{1}}(x_{1}-x_{2})^{j}A(x_{1},g^{-1}(x_{2})).$$
As $A(x_{1},g^{-1}(x_{2}))\in \Hom \left(W,W((x_{2}^{-1}))[[x_{1},x_{1}^{-1}]]\right)$, we get
$B_{j}(x)\in \E^{o}(W)$. It is clear that
(\ref{eagb-exp}) is equivalent to
\begin{eqnarray}
A(x_{1},x_{2})=\sum_{j=0}^{k-1} B_{j}(g(x_{2}))\Phi(g)^{-j}
\left(\frac{\partial}{\partial x_{2}}\right)^{j}x_{1}^{-1}\delta\left(\frac{g(x_{2})}{x_{1}}\right)
\end{eqnarray}
 under the condition that $B_{j}(x)\in \E^{o}(W)$ for $0\le j\le k-1$.

We next consider the general case. For $1\le i\le r$, set
$$p_{i}(x_{1},x_{2})=\frac{p(x_{1},x_{2})}{(x_{1}-g_{i}(x_{2}))^{k_{i}}}
\in \C[x_{1},x_{2}].$$
Note that $p_{i}(x_{1},x_{2})$ with $1\le i\le r$, viewed as polynomials in $x_{1}$ with coefficients in $\C(x_{2})$,
are relatively prime. Then there exist
$q_{i}(x_{1},x_{2})\in \C(x_{2})[x_{1}]$ for $1\le i\le r$ such that
\begin{eqnarray*}
1=p_{1}(x_{1},x_{2})q_{1}(x_{1},x_{2})+\cdots +p_{k}(x_{1},x_{2})q_{k}(x_{1},x_{2}).
\end{eqnarray*}
Recall the field embedding $\iota_{x,\infty}:\ \C(x)\rightarrow \C((x^{-1}))$. This gives an algebra embedding of $\C(x_{2})[x_{1}]$
into $\C((x_{2}^{-1}))[x_{1}]$.
Denote by $\bar{q}_{i}(x_{1},x_{2})$ the image of $q_{i}(x_{1},x_{2})$ in $\C((x_{2}^{-1}))[x_{1}]$.
We have
\begin{eqnarray*}
1=p_{1}(x_{1},x_{2})\bar{q}_{1}(x_{1},x_{2})+\cdots +p_{k}(x_{1},x_{2})\bar{q}_{k}(x_{1},x_{2})
\end{eqnarray*}
in $\C((x_{2}^{-1}))[x_{1}]$, so that
\begin{eqnarray}\label{e1-expansion}
A(x_{1},x_{2})=\sum_{i=1}^{k}p_{i}(x_{1},x_{2})\bar{q}_{i}(x_{1},x_{2})A(x_{1},x_{2}).
\end{eqnarray}
Note that
$$(x_{1}-g_{i}(x_{2}))^{k_{i}}p_{i}(x_{1},x_{2})\bar{q}_{i}(x_{1},x_{2})A(x_{1},x_{2})=0$$
and $p_{i}(x_{1},x_{2})\bar{q}_{i}(x_{1},x_{2})A(x_{1},x_{2})\in \Hom \left(W,W((x_{2}^{-1}))[[x_{1},x_{1}^{-1}]]\right).$
Then it follows immediately from (\ref{e1-expansion}) and the special case.
\end{proof}

The following is a variation of the notion of $\Gamma$-vertex algebra,
introduced in \cite{li-infinity} (cf. \cite{li-new})\footnote{The notion of
vertex $\Gamma$-algebra here with respect to the pair $(\Gamma,\phi)$ amounts
to the notion of $\Gamma$-vertex algebra in \cite{li-infinity} with respect to $(\Gamma,\phi^{-1})$.}:

\bd{dvertex-gamma-algebra}
{\em Let $\Gamma$ be a group. A {\em vertex $\Gamma$-algebra} is a vertex algebra $V$ equipped with
two group homomorphisms
$$L:\ \Gamma\rightarrow GL(V)\ \ \mbox{ and }\ \ \phi: \ \Gamma\rightarrow \C^{\times},$$
satisfying the condition that $L(g){\bf 1}={\bf 1}$ for $g\in \Gamma$ and
\begin{eqnarray}
L(g)Y(u,x)v=Y\left(L(g)u,\phi(g)x\right)L(g)v\ \ \ \mbox{ for }g\in \Gamma,\ u,v\in V.
\end{eqnarray}}
\ed

It is clear that for any group $H$ equipped with a homomorphism from $H$ to $\Gamma$,
a vertex $\Gamma$-algebra $V$ is naturally a vertex $H$-algebra.
Define a {\em homomorphism of vertex $\Gamma$-algebras} from $U$ to $V$
to be a homomorphism $\theta$ of vertex algebras such that
$$\theta (L(g)u)=L(g)\theta(u)\ \ \ \mbox{ for }g\in \Gamma,\ u\in U.$$

\br{rZgradedva}
{\em Let $V$ be a $\Z$-graded vertex algebra in the sense that $V$ is a vertex algebra
equipped with a $\Z$-grading $V=\oplus_{n\in \Z}V_{(n)}$ such that
$$u_{k}v\in V_{(m+n-k-1)}\ \ \ \mbox{ for }u\in V_{(m)},\ v\in V_{(n)},\ m,n,k\in \Z.$$
Let $L(0)$ denote the linear operator on $V$, defined by
$$L(0)v=nv\ \ \ \mbox{ for }v\in V_{(n)},\ n\in \Z.$$
Let $\Gamma$ be any automorphism group of $V$, which preserves the
$\Z$-grading, and let $\phi: \Gamma\rightarrow \C^{\times}$ be any
linear character. For $g\in \Gamma$, set $L(g)=\phi(g)^{L(0)}g$.
Then it can be readily seen that $V$ becomes a vertex
$\Gamma$-algebra.} \er

The following is a modification of the same named notion
introduced in \cite{li-infinity}\footnote{The only change from the original definition is in (\ref{ecovariance}),
to incorporate the change of the right action $R$ of $G$ on $\E^{o}(W)$
to the left action $L$.}:

\bd{dquasi-module-vgalgebra}
{\em Let $V$ be a vertex $\Gamma$-algebra. A {\em quasi $V$-module at infinity} is
a quasi module at infinity $(W,Y_{W})$ for $V$ viewed as a vertex algebra, equipped with
a group homomorphism
$$\Psi: \ \Gamma\rightarrow G,$$
satisfying the condition that $\phi=\Phi \circ \Psi$,
\begin{eqnarray}\label{ecovariance}
Y_{W}(L(g)v,x)=Y_{W}(v,\Psi(g)^{-1}(x))\ \ \ \mbox{ for }g\in
\Gamma,\ v\in V,
\end{eqnarray}
and $\{ Y_{W}(v,x)\ |\ v\in V\}$ is $\Psi(\Gamma)$-local.}
\ed

\bl{lweaker-compatibility}
In Definition \ref{dquasi-module-vgalgebra}, the condition that
$\{ Y_{W}(v,x)\ |\ v\in V\}$ is $\Psi(\Gamma)$-local can be replaced with a weaker condition that
for any $u,v\in V$, there exists a product $q(x_{1},x_{2})$ of linear polynomials $x_{1}-g(x_{2})$ with $g\in \Psi(\Gamma)$ such that
\begin{eqnarray}\label{eqpgamma-comp}
q(x_{1},x_{2})Y_{W}(u,x_{1})Y_{W}(v,x_{2})\in \Hom \left(W,W((x_{1}^{-1},x_{2}^{-1}))\right).
\end{eqnarray}
\el

\begin{proof} Let $u,v\in V$. Assume that $q(x_{1},x_{2})$ is a product of linear polynomials $x_{1}-g(x_{2})$ with $g\in \Psi(\Gamma)$
such that (\ref{eqpgamma-comp}) holds. {}From Proposition \ref{pjacobi}, there exists a nonzero polynomial $p(x_{1},x_{2})$ such that
$$p(x_{1},x_{2})Y_{W}(u,x_{1})Y_{W}(v,x_{2})=p(x_{1},x_{2})Y_{W}(v,x_{2})Y_{W}(u,x_{1}).$$
Then
\begin{eqnarray*}
 p(x_{1},x_{2})\left(q(x_{1},x_{2})Y_{W}(u,x_{1})Y_{W}(v,x_{2})-q(x_{1},x_{2})Y_{W}(v,x_{2})Y_{W}(u,x_{1})\right)=0.
\end{eqnarray*}
{}From definition, $q(x_{1},x_{2})Y_{W}(v,x_{2})Y_{W}(u,x_{1})$ lies in
 $\Hom \left(W,W((x_{2}^{-1}))((x_{1}^{-1}))\right)$.
As (\ref{eqpgamma-comp}) holds, $q(x_{1},x_{2})Y_{W}(u,x_{1})Y_{W}(v,x_{2})$ also lies in
$\Hom \left(W,W((x_{2}^{-1}))((x_{1}^{-1}))\right)$.
Consequently, (multiplying by $\iota_{x_{2}^{-1},x_{1}^{-1}}(1/p(x_{1},x_{2}))$) we get
$$q(x_{1},x_{2})Y_{W}(u,x_{1})Y_{W}(v,x_{2})=q(x_{1},x_{2})Y_{W}(v,x_{2})Y_{W}(u,x_{1}).$$
This proves that $\{ Y_{W}(v,x)\ |\ v\in V\}$ is $\Psi(\Gamma)$-local.
\end{proof}

Recall the (left) group action of $G$ on $\E^{o}(W)$. We have
(cf. \cite{li-infinity}, Theorem 5.10)\footnote{
This is the corrected version of Theorem 5.10 in \cite{li-infinity} with the right action $R$ of $\Gamma$
replaced with the left action $L$.}:

\bt{tabstract} Let $W$ be a vector space, let $\Gamma$ be a subgroup
of $G$, and let $U$ be a $\Gamma$-local subset of $\E^{o}(W)$. Then
$\Gamma\cdot U$ is $\Gamma$-local and the vertex algebra
$\<\Gamma\cdot U\>$ generated by $\Gamma \cdot U$ is a vertex
$\Gamma$-algebra
 where  $\phi=\Phi$ and $L$ is given by $L(g)=L_{g}$ for $g\in \Gamma$, that is
 $$L(g)a(x)=a(g^{-1}(x))\ \ \mbox{ for }g\in \Gamma,\ a(x)\in\<\Gamma\cdot
 U\>.$$ Furthermore, $W$ is a quasi module at infinity
 for $\<\Gamma\cdot U\>$ with $Y_{W}(a(x),z)=a(z)$ for
$a(x)\in\<\Gamma\cdot U\>$ and with $\Psi=1$. \et

\begin{proof} As $L_{g}=R_{g^{-1}}$ on $\E^{o}(W)$
for $g\in G$, Lemma 3.13 of \cite{li-infinity} asserts that
$$L_{g}Y_{\E^{o}}(a(x),x_{0})b(x)=Y_{\E^{o}}(L_{g}a(x),\Phi(g)x_{0})L_{g}b(x)$$
for $g\in \Gamma,\ a(x),b(x)\in \<\Gamma\cdot U\>$.
This together with Theorem \ref{trecall-main} confirms the first assertion on the vertex $\Gamma$-algebra structure.
As for the structure of a quasi module at infinity, we have
$$Y_{W}(L_{g}a(x),x_{0})=Y_{W}(a(g^{-1}(x)),x_{0})=a(g^{-1}(x_{0}))=Y_{W}(a(x),g^{-1}(x_{0}))$$
for $g\in \Gamma,\ a(x)\in \<\Gamma\cdot U\>$. On the other hand, it follows from Proposition 3.14 of \cite{li-infinity} and induction
that $\< \Gamma\cdot U\>$ is $\Gamma$-local. This confirms the second assertion that $W$ is a quasi module at infinity
for $\<\Gamma \cdot U\>$ viewed as a vertex $\Gamma$-algebra with $Y_{W}(a(x),x_{0})=a(x_{0})$ for $a(x)\in \<\Gamma\cdot U\>$.
\end{proof}

The following refinement of Proposition \ref{prepresentation} is straightforward:

\bp{prepresentation-refine}
Let $V$ be a vertex $\Gamma$-algebra, let $\Psi: \Gamma \rightarrow G$ be a group homomorphism such that
$\Phi\circ \Psi=\phi$,
 and let $W$ be a vector space equipped with a linear map
$Y_{W}(\cdot,x): V\rightarrow \Hom (W,W((x^{-1})))$ with $Y_{W}({\bf 1},x)=1_{W}$.
 Set
$$V_{W}=\{ Y_{W}(v,x)\ |\ v\in V\}\subset \E^{o}(W).$$
Then
$(W,Y_{W})$ is a quasi $V$-module at infinity if and only if $V_{W}$ is $\Psi(\Gamma)$-local,
$(V_{W},Y_{\E^{o}},1_{W})$ carries the structure of a vertex $\Psi(\Gamma)$-algebra,
 and $Y_{W}: V\rightarrow V_{W}$ is a homomorphism of vertex $\Gamma$-algebras, where $V_{W}$ is viewed as a vertex $\Gamma$-algebra
 through the homomorphism $\Psi: \Gamma\rightarrow \Psi(\Gamma)\subset G$.
\ep

As the main result of this section, we have
the following analog of the twisted vertex operator commutator formula (see \cite{flm}):

\bt{tconverse-general}
Let $V$ be a vertex $\Gamma$-algebra with group homomorphisms $L: \Gamma\rightarrow GL(V);\ \phi: \Gamma\rightarrow \C^{\times}$,
let $(W,Y_{W})$ be a quasi $V$-module at infinity with group homomorphism $\Psi:\Gamma\rightarrow G$, and
let $u,v\in V$. Then there exist finitely many $\sigma_{1},\dots,\sigma_{r}\in \Gamma$  such that
$\Psi(\sigma_{i})$ $(i=1,\dots,r)$ are distinct in $G$ and
\begin{eqnarray}
&&[Y_{W}(u,x_{1}),Y_{W}(v,x_{2})]\\
&=&-\sum_{i=1}^{r}\sum_{j\in \N}\phi(\sigma_{i})^{-1} Y_{W}\left((L(\sigma_{i})u)_{j}v,x_{2}\right)\frac{1}{j!}\left(\frac{\partial}{\partial x_{2}}\right)^{j}x_{1}^{-1}\delta\left(\frac{\Psi(\sigma_{i})^{-1}(x_{2})}{x_{1}}\right).\nonumber
\end{eqnarray}
Furthermore, if $\Psi: \Gamma\rightarrow G$ is one-to-one, then
\begin{eqnarray}
&&[Y_{W}(u,x_{1}),Y_{W}(v,x_{2})]\\
&=&-\sum_{\sigma\in \Gamma,j\in \N} \phi(\sigma)^{-1}Y_{W}\left((L(\sigma)u)_{j}v,x_{2}\right)\frac{1}{j!}\left(\frac{\partial}{\partial x_{2}}\right)^{j}x_{1}^{-1}\delta\left(\frac{\Psi(\sigma)^{-1}(x_{2})}{x_{1}}\right),\nonumber
\end{eqnarray}
which is a finite sum.
\et

\begin{proof}
{}From definition, there exist distinct $g_{1}(x),\dots, g_{r}(x)\in \Psi(\Gamma)\subset G$ and positive integers
$k_{1},\dots,k_{r}$ such that
\begin{eqnarray}\label{epoly-comm}
\left(\prod_{i=1}^{r}(x_{1}-g_{i}(x_{2}))^{k_{i}}\right)[Y_{W}(u,x_{1}),Y_{W}(v,x_{2})]=0.
\end{eqnarray}
In view of Lemma \ref{lformal-variable}, we have
\begin{eqnarray}\label{eywuv-aij}
\mbox{}\ \ \ \ \ \ \ \ [Y_{W}(u,x_{1}),Y_{W}(v,x_{2})]
=\sum_{i=1}^{r}\sum_{j=0}^{k_{i}-1}A_{ij}(x_{2})\frac{1}{j!}\left(\frac{\partial}{\partial x_{2}}\right)^{j}
x_{1}^{-1}\delta\left(\frac{g_{i}(x_{2})}{x_{1}}\right)
\end{eqnarray}
for some $A_{ij}(x)\in \E^{o}(W)$. By Proposition \ref{pextension},
we get
$$A_{ij}(x)=-\Phi(g_{i}) Y_{W}(u,g_{i}(x))_{j}Y_{W}(v,x),$$
where $Y_{W}(u,g_{i}(x))$ and $Y_{W}(v,x)$ are viewed as elements of
$\E^{o}(W)$. Let $\sigma_{1},\dots,\sigma_{r}\in \Gamma$ such that
$\Psi(\sigma_{i})=g_{i}^{-1}(x)$ for $1\le i\le r$. As
$Y_{W}(L(\sigma_{i})u,x)=Y_{W}(u,g_{i}(x))$ and $\Phi(g_{i})=\Phi(\Psi(\sigma_{i}^{-1}))=\phi(\sigma_{i})^{-1}$, we have
$$A_{ij}(x)=-\phi(\sigma_{i})^{-1} Y_{W}(L(\sigma_{i})u,x)_{j}Y_{W}(v,x).$$
 Recall from Proposition \ref{prepresentation} that
\begin{eqnarray}\label{ehom-property}
Y_{\E^{o}}(Y_{W}(a,x),x_{0})Y_{W}(b,x)=Y_{W}(Y(a,x_{0})b,x)\ \ \ \mbox{ for }a,b\in V.
\end{eqnarray}
Using this we obtain
\begin{eqnarray*}
&&[Y_{W}(u,x_{1}),Y_{W}(v,x_{2})]\nonumber\\
&=&-\Res_{x_{0}}\sum_{i=1}^{r}\Phi(g_{i})Y_{\E^{o}}\left( Y_{W}(L(\sigma_{i})u,x_{2}),x_{0}\right)Y_{W}(v,x_{2})
e^{x_{0}\frac{\partial}{\partial x_{2}}}
x_{1}^{-1}\delta\left(\frac{g_{i}(x_{2})}{x_{1}}\right)\nonumber\\
&=&-\Res_{x_{0}}\sum_{i=1}^{r}\phi(\sigma_{i})^{-1}Y_{W}\left(Y(L(\sigma_{i})u,\Phi(g_{i})^{-1}x_{0})v,x_{2}\right)
e^{x_{0}\frac{\partial}{\partial x_{2}}}x_{1}^{-1}\delta\left(\frac{g_{i}(x_{2})}{x_{1}}\right)\nonumber\\
&=&-\Res_{x_{0}}\sum_{i=1}^{r}\phi(\sigma_{i})^{-1}Y_{W}\left(Y(L(\sigma_{i})u,\phi(\sigma_{i})x_{0})v,x_{2}\right)
e^{x_{0}\frac{\partial}{\partial x_{2}}}x_{1}^{-1}\delta\left(\frac{g_{i}(x_{2})}{x_{1}}\right),\ \ \ \
\end{eqnarray*}
as desired.

Suppose $\Psi$ is one-to-one.
Let $\sigma\in \Gamma$ with $\sigma\ne \sigma_{i}$ for $1\le i\le r$. Set $g(x)=\Psi(\sigma^{-1})\in G$.
Then $g(x)\ne g_{i}(x)$ for $1\le i\le r$. Combining (\ref{eywuv-aij}) with
 Proposition \ref{pextension}, we have
 $$Y_{W}(u,g(x))_{j}Y_{W}(v,x)=0\ \ \mbox{ for }j\ge 0,$$
i.e.,
\begin{equation*}
Y_{\E^o}\left(Y_{W}(u,g(x)),x_0\right)Y_{W}(v,x)\in \E^{o}(W)[[x_0]].
\end{equation*}
Using Proposition \ref{prepresentation} we have
\begin{eqnarray*}
Y_{W}\left(Y(L(\sigma)u,x_{0})v,x\right)
&=&Y_{\E^{o}}\left(Y_{W}(L(\sigma)u,x),x_{0}\right)Y_{W}(v,x)\\
&=&Y_{\E^{o}}\left(Y_{W}(u,g(x)),x_{0}\right)Y_{W}(v,x).
\end{eqnarray*}
Consequently, $Y_{W}\left(Y(L(\sigma)u,x_{0})v,x\right)$ involves
only nonnegative powers of $x_{0}$. Thus
$$\Res_{x_{0}}Y_{W}(Y(L(\sigma)u,x_{0})v,x_{2})e^{x_{0}\frac{\partial}{\partial x_{2}}}
x_{1}^{-1}\delta\left(\frac{g(x_{2})}{x_{1}}\right) =0.$$ Then the
second assertion follows immediately.
\end{proof}

The following technical result is a modification of Lemma 5.9 in \cite{li-infinity}
with a slightly different proof:

\bl{lva-gamma-quasi-module} Let $V$ be a vertex $\Gamma$-algebra,
let $\Psi: \Gamma\rightarrow G$ be a group homomorphism with
$\Phi\circ \Psi=\phi$, and let $(W,Y_{W})$ be a quasi module at infinity
for $V$ viewed as a vertex algebra. Assume that $U_{W}=\{
Y_{W}(u,x)\;|\; u\in U\}$ is $\Psi(\Gamma)$-local and
$$Y_{W}(L(g)u,x)=Y_{W}(u,\Psi(g)^{-1}(x))\ \ \ \mbox{ for }g\in \Gamma,\; u\in U,$$
where $U$ is a $\Gamma$-submodule which generates $V$ as a vertex algebra.
Then $(W,Y_{W})$ is a quasi $V$-module at infinity. \el

\begin{proof} {}From assumption,
$U_{W}$ is a $\Psi(\Gamma)$-local subset of $\E^{o}(W)$. For $g\in \Gamma,\ u\in U$, we have
\begin{eqnarray}\label{ecovariance-need}
Y_{W}(L(g)u,x)=Y_{W}(u,\Psi(g)^{-1}(x))=L_{\Psi(g)}\left(Y_{W}(u,x)\right).
\end{eqnarray}
It follows that $U_{W}$ is stable under the subgroup $\Psi(\Gamma)$ of $G$.
By Theorem \ref{tabstract}, we have a vertex $\Psi(\Gamma)$-algebra $\<U_{W}\>$ with $W$ as a quasi module at infinity.
Since $(W,Y_{W})$ is a module at infinity for $V$ viewed as a vertex algebra,
in view of Proposition \ref{prepresentation},  $Y_{W}$ is a homomorphism of vertex algebras from $V$ to $\<U_{W}\>$.

Suppose that
$$Y_{W}(L(g)u,x)=Y_{W}(u,\Psi(g)^{-1}(x))\ \mbox{ and }
Y_{W}(L(g)v,x)=Y_{W}(v,\Psi(g)^{-1}(x))$$ for some $g\in \Gamma,\;
u,v\in V$. Then we have
\begin{eqnarray*}
& &Y_{W}(L(g)Y(u,x_{0})v,x)\\
&=&Y_{W}\left(Y(L(g)u,\phi(g)x_{0})L(g)v,x\right)\\
&=&Y_{\E^{o}}\left(Y_{W}(L(g)u,x),\phi(g)x_{0}\right)Y_{W}(L(g)v,x)\\
&=&Y_{\E^{o}}\left(Y_{W}(u,\Psi(g)^{-1}(x)),\phi(g)x_{0}\right)Y_{W}(v,\Psi(g)^{-1}(x))\\
&=&L_{\Psi(g)}Y_{\E^{o}}\left(Y_{W}(u,x),x_{0}\right)Y_{W}(v,x)\\
&=&L_{\Psi(g)}Y_{W}\left(Y(u,x_{0})v,x\right)\\
&=&Y_{W}(Y(u,x_{0})v,\Psi(g)^{-1}(x)).
\end{eqnarray*}
As $U$ generates $V$ as a vertex algebra, it follows from (\ref{ecovariance-need}) that
$$Y_{W}(L(g)v,x)=L_{\Psi(g)}\left(Y_{W}(v,x)\right)\ \ \ \mbox{ for all }g\in \Gamma,\ v\in V.$$
That is, $Y_{W}$ is a homomorphism of vertex $\Gamma$-algebras.
By Proposition \ref{prepresentation-refine}, $(W,Y_{W})$ is a quasi $V$-module at infinity.
\end{proof}

\section{Lie algebra $\hat{\g}(\infty)[\Gamma]$ and vertex algebra $V_{\hat{\g}}(\ell,0)$}
In this section, we recall from \cite{li-infinity} the Lie algebra $\hat{\g}(\infty)[\Gamma]$
and the main results on the relation between $\hat{\g}(\infty)[\Gamma]$ and vertex algebra $V_{\hat{\g}}(\ell,0)$,
including Theorem 5.14, and as our main result we establish the converse of this theorem.

Let $\g$ be a (possibly infinite-dimensional) Lie algebra equipped with a
non-degenerate symmetric invariant bilinear form $\<\cdot,\cdot\>$.
Associated to the pair $(\g,\<\cdot,\cdot\>)$, one has an
(untwisted) affine Lie algebra
$$\hat{\g}=\g\otimes \C[t,t^{-1}]\oplus \C {\bf k},$$
where ${\bf k}$ is central and
$$[a\otimes t^{m},b\otimes
t^{n}]=[a,b]\otimes t^{m+n}+m\delta_{m+n,0}\<a,b\>{\bf k}$$ for
$a,b\in \g,\; m,n\in \Z$. Defining $\deg (\g\otimes t^{m})=-m$ for
$m\in \Z$ and $\deg {\bf k}=0$ makes $\hat{\g}$ a $\Z$-graded Lie
algebra. For $a\in \g$, form a generating function
 $$a(x)=\sum_{n\in \Z}(a\otimes t^{n})x^{-n-1}.$$

Let $\ell$ be a complex number. Denote by $\C_{\ell}$ the
$1$-dimensional $(\g\otimes \C[t]\oplus \C {\bf k})$-module with
$\g\otimes \C[t]$ acting trivially and with ${\bf k}$ acting as
scalar $\ell$. Form the induced $\hat{\g}$-module
\begin{eqnarray}
V_{\hat{\g}}(\ell,0)=U(\hat{\g})\otimes_{U(\g\otimes \C[t]
\oplus \C {\bf k})} \C_{\ell}.
\end{eqnarray}
 Set ${\bf 1}=1\otimes 1$ and then identify $\g$ as
a subspace of $V_{\hat{\g}}(\ell,0)$ through the linear map
$a\rightarrow a(-1){\bf 1}$. It is well known (cf. \cite{fz})
that there exists a unique vertex-algebra structure on
$V_{\hat{\g}}(\ell,0)$ with ${\bf 1}$ as the vacuum vector and with
$Y(a,x)=a(x)$ for $a\in \g$. Defining $\deg {\bf 1}=0$ makes
$V_{\hat{\g}}(\ell,0)$ a $\Z$-graded $\hat{\g}$-module and the
vertex algebra $V_{\hat{\g}}(\ell,0)$ equipped with this
$\Z$-grading is a $\Z$-graded vertex algebra.

Let $\Gamma$ be a subgroup of $\Aut(\g,\<\cdot,\cdot\>)$, consisting of automorphisms of $\g$ that preserve $\<\cdot,\cdot\>$.
Each $g\in \Gamma$ lifts canonically to an
automorphism of the $\Z$-graded Lie algebra $\hat{\g}$, and then to an automorphism of
the $\Z$-graded vertex algebra $V_{\hat{\g}}(\ell,0)$. In this way,
$\Gamma$ acts on vertex algebra $V_{\hat{\g}}(\ell,0)$ by
automorphisms that preserve the $\Z$-grading. Let $\phi:
\Gamma\rightarrow \C^{\times}$ be any group homomorphism. For $g\in
\Gamma$, set
$$L(g)=\phi(g)^{L(0)}g\in GL(V_{\hat{\g}}(\ell,0)),$$
where $L(0)$ denotes the linear operator on $V_{\hat{\g}}(\ell,0)$
defined by $L(0)v=nv$ for $v\in V_{\hat{\g}}(\ell,0)_{(n)}$ with $n\in \Z$. This defines a
vertex $\Gamma$-algebra structure on $V_{\hat{\g}}(\ell,0)$.

Consider the following completion of affine Lie
algebra $\hat{\g}$:
\begin{eqnarray}
\hat{\g}(\infty)=\g\otimes \C((t^{-1}))\oplus \C {\bf k},
\end{eqnarray}
where
\begin{eqnarray}
[a\otimes p(t),b\otimes q(t)]=[a,b]\otimes
p(t)q(t)+\Res_{t}p'(t)q(t)\<a,b\>{\bf k}
\end{eqnarray}
for $a,b\in \g,\; p(t),q(t)\in \C((t^{-1}))$.

The following is a construction of a family of new Lie algebras by using
Lie algebra $\hat{\g}(\infty)$
(cf. \cite{li-infinity}, Proposition 5.12; \cite{gkk})\footnote{This is the corrected version of Proposition 5.12 in \cite{li-infinity}.
In the original proof, an action of $\Gamma$ on $\hat{\g}(\infty)$ was defined by
$$g(a\otimes p(t)+\lambda {\bf k})
=ga\otimes p(g(t))+\lambda {\bf k}$$
for $g\in\Gamma,\; a\in \g,\; p(t)\in \C((t^{-1})),\; \lambda\in \C$, which is not a left action
if $\Psi(\Gamma)$ is not abelian.}:

\bp{ptwisted-affine-comp}
Let $\g$ be a Lie algebra equipped with a
non-degenerate symmetric invariant bilinear form
 $\<\cdot,\cdot\>$ and let $\Gamma$ be an automorphism group of $(\g,\<\cdot,\cdot\>)$,
satisfying the condition that for any $u,v\in \g$,
$$[gu,v]=0\;\;\mbox{ and }\;\; \<gu,v\>=0\;\;\;\mbox{ for all but
finitely many }g\in \Gamma.$$ Let $\Psi: \Gamma\rightarrow G$ be a
group homomorphism and set $g(x)=\Psi(g)(x)\in G$ for $g\in \Gamma$.
Define a new bilinear multiplicative operation
$[\cdot,\cdot]_{\Gamma}$ on vector space
$\hat{\g}(\infty)=\g\otimes \C((t^{-1}))\oplus \C {\bf k}$ by
\begin{eqnarray*}
& &[a\otimes p(t),{\bf k}]_{\Gamma}=0=[{\bf k},a\otimes p(t)]_{\Gamma}, \\
& &[a\otimes p(t),b\otimes q(t)]_{\Gamma}=\sum_{g\in
\Gamma}[ga,b]\otimes p(g^{-1}(t))q(t)+
\Res_{t}q(t)\frac{d}{dt}p(g^{-1}(t))\<ga,b\>{\bf k}\ \ \ \ \ \
\end{eqnarray*}
for $a,b\in \g,\; p(t),q(t)\in \C((t^{-1}))$. Then the subspace,
linearly spanned by vectors $$ga\otimes p(t)-a\otimes p(g(t))$$
for $g\in \Gamma,\;a\in \g,\;p(t)\in\C((t^{-1}))$, is a two-sided
ideal of the non-associative algebra, and the quotient algebra, which
we denote by $\hat{\g}(\infty)[\Gamma]$, is a Lie algebra. \ep

\begin{proof}
For $g\in\Gamma,\; a\in \g,\; p(t)\in \C((t^{-1})),\; \lambda\in \C$,
define
$$g(a\otimes p(t)+\lambda {\bf k})
=ga\otimes p(g^{-1}(t))+\lambda {\bf k}.$$
It is straightforward to show that $\Gamma$ acts on $\hat{\g}(\infty)$ by
automorphisms, satisfying the condition that for any $u,v\in \hat{\g}(\infty)$,
$$[gu,v]=0\ \ \ \mbox{ for all but finitely many }g\in \Gamma.$$
Then it follows
immediately from \cite{li-gamodule} (Lemma 4.1).
\end{proof}

Let
$$\pi:\ \hat{\g}(\infty)\rightarrow \hat{\g}(\infty)[\Gamma] $$
denote the natural linear map. For $a\in \g$,  set
\begin{eqnarray}
a_{\Gamma}(x)=\sum_{n\in \Z}\pi (a\otimes t^{n}) x^{-n-1}
\in \left(\hat{\g}(\infty)[\Gamma]\right)[[x,x^{-1}]].
\end{eqnarray}
Define a linear character $\phi:\Gamma \rightarrow \C^{\times}$ by
$$\phi(g)=\frac{d}{dx}\Psi(g)(x)\ \ \mbox{ for }g\in \Gamma.$$

We say that a $\hat{\g}(\infty)[\Gamma]$-module $W$ is of {\em level} $\ell\in \C$
if ${\bf k}$ acts on $W$ as scalar $\ell$. The following is a modification of Lemma 5.13 in \cite{li-infinity}:

\bl{lgenerating-relation} Let $W$ be a vector space and let $\ell\in \C$.
Then a $\hat{\g}(\infty)[\Gamma]$-module structure of level $\ell$ on $W$
amounts to a linear map
$$\theta:\ \g \rightarrow \E^{o}(W);\ \ \ a\mapsto a_{W}(x),$$
satisfying the conditions that
\begin{eqnarray}\label{e5.15}
(ga)_{W}(x)=\phi(g)^{-1}a_{W}(g^{-1}(x))
\end{eqnarray}
for $g\in \Gamma,\; a\in \g$ and that
\begin{eqnarray}\label{ecomm-liealgebra}
&&[a_{W}(x_{1}),b_{W}(x_{2})]\\
&=&\sum_{g\in
\Gamma}[ga,b]_{W}(x_{2})
x_{1}^{-1}\delta\left(\frac{g^{-1}(x_{2})}{x_{1}}\right) +\ell\<ga,b\>\frac{\partial}{\partial
x_{2}}x_{1}^{-1}\delta\left(\frac{g^{-1}(x_{2})}{x_{1}}\right)\nonumber
\end{eqnarray}
for $a,b\in \g$. \el

\begin{proof}
Let $g(x)=\alpha x+\beta\in G$ with $\alpha\in \C^{\times},\;
\beta\in \C$. Then $g^{-1}(x)=\alpha^{-1}(x-\beta)$ and
\begin{eqnarray*}
x^{-1}\delta\left(\frac{g(t)}{x}\right)
=x^{-1}\delta\left(\frac{\alpha t+\beta}{x}\right) =\alpha^{-1}
t^{-1}\delta\left(\frac{x-\beta}{\alpha t}\right) =\alpha^{-1}
t^{-1}\delta\left(\frac{g^{-1}(x)}{t}\right).
\end{eqnarray*}
Using this we get
\begin{eqnarray*}
(ga)_{\Gamma}(x)&=&\sum_{n\in\Z}\pi(ga\otimes t^{n})x^{-n-1} =
\sum_{n\in\Z}\pi(a\otimes
g(t)^{n})x^{-n-1}\\
&=&\pi\left(a\otimes x^{-1}\delta\left(\frac{g(t)}{x}\right)\right)\\
&=&\phi(g)^{-1} \pi\left(a\otimes
t^{-1}\delta\left(\frac{g^{-1}(x)}{t}\right)\right)\\
&=&\phi(g)^{-1}a_{\Gamma}(g^{-1}(x)),
\end{eqnarray*}
proving (\ref{e5.15}). As for (\ref{ecomm-liealgebra}), notice that
\begin{eqnarray*}
&&\sum_{m,n\in \Z}g^{-1}(t)^{m}t^{n}x_{1}^{-m-1}x_{2}^{-n-1}
=x_{1}^{-1}\delta\left(\frac{g^{-1}(t)}{x_{1}}\right)
x_{2}^{-1}\delta\left(\frac{t}{x_{2}}\right)\\
&&=x_{1}^{-1}\delta\left(\frac{g^{-1}(x_{2})}{x_{1}}\right)
x_{2}^{-1}\delta\left(\frac{t}{x_{2}}\right)
\end{eqnarray*}
 and
\begin{eqnarray*}
&&\Res_{t}\sum_{m,n\in
\Z}t^{n}\frac{d}{dt}g^{-1}(t)^{m}
x_{1}^{-m-1}x_{2}^{-n-1}\\
&=&-\Res_{t}\sum_{m,n\in \Z}nt^{n-1}g^{-1}(t)^{m}x_{1}^{-m-1}x_{2}^{-n-1}\\
&=&\Res_{t}\frac{\partial}{\partial x_{2}}
x_{1}^{-1}\delta\left(\frac{g^{-1}(t)}{x_{1}}\right)
t^{-1}\delta\left(\frac{x_{2}}{t}\right)\\
&=&\frac{\partial}{\partial x_{2}}
x_{1}^{-1}\delta\left(\frac{g^{-1}(x_{2})}{x_{1}}\right).
\end{eqnarray*}
Then (\ref{ecomm-liealgebra}) follows from the construction of $\hat{\g}(\infty)[\Gamma]$.
\end{proof}

We have (cf. \cite{li-infinity}, Theorem 5.14)\footnote{This is a corrected version of
Theorem 5.14 in \cite{li-infinity} with a complete proof.}:

\bt{t5.14}
Let $\ell\in \C$ and let $W$ be any $\hat{\g}(\infty)[\Gamma]$-module of level $\ell$
such that $a_{\Gamma}(x)\in \E^{o}(W)$ for $a\in \g$. Then there exists a unique structure of
a quasi module at infinity on $W$ for the vertex $\Gamma$-algebra $V_{\hat{\g^{o}}}(-\ell,0)$ with
$Y_{W}(a,x)=a_{\Gamma}(x)$ for $a\in \g$, where $\g^{o}$ denotes the opposite Lie algebra of $\g$.
\et

\begin{proof} Set
$$U=\{ a_{\Gamma}(x)\;|\; a\in \g\}\subset \E^{o}(W).$$
For $a,b\in \g$, there exist $g_{1},\dots,g_{r}\in \Gamma$ such that
$[ga,b]=0$ and $\<ga,b\>=0$ for $g\notin \{g_{1},\dots,g_{r}\}$.
It follows from (\ref{ecomm-liealgebra}) that
$$(x_{1}-g_{1}^{-1}(x_{2}))^{2}\cdots
(x_{1}-g_{r}^{-1}(x_{2}))^{2}[a_{\Gamma}(x_{1}),b_{\Gamma}(x_{2})]=0.$$
Thus $U$ is a $\Gamma$-local subspace of $\E^{o}(W)$. {}From
(\ref{e5.15}), $\Gamma\cdot U=U$. By Theorem \ref{tabstract}, $U$
generates a vertex $\Gamma$-algebra $\<U\>$ with $W$ as a quasi
module-at-infinity where $Y_{W}(\alpha(x),x_{0})=\alpha(x_{0})$ for $\alpha(x)\in \<U\>$.
Combining (\ref{ecomm-liealgebra}) with Lemma \ref{lextension}
we get
\begin{eqnarray*}
a_{\Gamma}(x)_{0}b_{\Gamma}(x)=-[a,b]_{\Gamma}(x),\;\;
a_{\Gamma}(x)_{1}b_{\Gamma}(x)=-\ell \<a,b\>1_{W}, \;\mbox{ and
}\;a_{\Gamma}(x)_{n}b_{\Gamma}(x)=0
\end{eqnarray*}
for $n\ge 2$. In view of the universal property of $V_{\hat{\g^{o}}}(-\ell,0)$ (cf.
\cite{li-gamodule}), there exists a vertex-algebra homomorphism from
$V_{\hat{\g^{o}}}(-\ell,0)$ to $\<U\>$, sending $a$ to
$a_{\Gamma}(x)$ for $a\in \g$. Consequently, $W$ is a quasi
module-at-infinity for $V_{\hat{\g^{o}}}(-\ell,0)$ viewed as a
vertex algebra. Furthermore, for $g\in \Gamma,\; a\in \g$ we have
$$Y_{W}(L(g)a,x)=Y_{W}(\phi(g)^{L(0)}ga,x)=\phi(g)(ga)_{\Gamma}(x)
=a_{\Gamma}(g^{-1}(x))=Y_{W}(a,g^{-1}(x)).$$ As $\g$ generates
$V_{\widehat{\g^{o}}}(-\ell,0)$ as a vertex algebra, it follows from
Lemma \ref{lva-gamma-quasi-module} that $W$ is a quasi
module-at-infinity for $V_{\hat{\g^{o}}}(-\ell,0)$ viewed as a
vertex $\Gamma$-algebra.
\end{proof}

Furthermore, we have the following converse of Theorem \ref{t5.14}:

\bt{tconverse}
Let $(W,Y_{W})$ be any quasi module at infinity for the vertex $\Gamma$-algebra $V_{\hat{\g^{o}}}(-\ell,0)$ such that
the associated homomorphism $\Psi: \Gamma\rightarrow G$ is one-to-one.
Then $W$ is a module for Lie algebra $\hat{\g}(\infty)[\Gamma]$ of level $\ell$ with $a_{\Gamma}(x)=Y_{W}(a,x)$ for $a\in \g$.
\et

\begin{proof} Note that for $u,v\in \g\subset V_{\hat{\g^{o}}}(-\ell,0)$, we have
$$u_{0}v=-[u,v],\ \ u_{1}v=-\ell \<u,v\>{\bf 1},\mbox{ and } \ u_{i}v=0\ \ \mbox{ for }i\ge 2.$$
Let $a,b\in \g$. Using Theorem \ref{tconverse-general} and the facts above we get
\begin{eqnarray*}
&&[Y_{W}(a,x_{1}),Y_{W}(b,x_{2})]\\
&=&-\sum_{g\in \Gamma}Y_{W}([ga,b],x_{2})x_{1}^{-1}\delta\left(\frac{g^{-1}(x_{2})}{x_{1}}\right)
+\ell \<ga,b\>\frac{\partial}{\partial x_{2}}x_{1}^{-1}\delta\left(\frac{g^{-1}(x_{2})}{x_{1}}\right),
\end{eqnarray*}
noticing that $L(g)=\phi(g)^{L(0)}g$ and $L(0)a=a$.
On the other hand, we have
$$Y_{W}(ga,x)=\phi(g)^{-1}Y_{W}(L(g)a,x)=\phi(g)^{-1}Y_{W}(a,g^{-1}(x))\ \ \mbox{ for }g\in \Gamma,\ a\in \g.$$
It then follows from Lemma \ref{lgenerating-relation}
that $W$ is a module for $\hat{\g}(\infty)[\Gamma]$ of level $\ell$ with $a_{\Gamma}(x)=Y_{W}(a,x)$ for $a\in \g$.
\end{proof}

\vspace{1cm}

\begin{center}
{\large \bf Acknowledgements}
\end{center}

The first named author (H.L) gratefully acknowledges the partial financial support from  NSA grant
H98230-11-1-0161 and China NSF grant 11128103.

\end{document}